\documentclass[12pt]{article}

\usepackage{amsmath,amssymb,amsfonts,theorem,makeidx,latexsym,epsfig}

\usepackage{color}

\newtheorem{defn}{Definition}[section]

\newtheorem{lemma}[defn]{Lemma}

{\theorembodyfont{\rmfamily}

\newtheorem{ex}[defn]{Example}}

\newtheorem{thm}[defn]{Theorem}

\newtheorem{prop}[defn]{Proposition}

\newtheorem{cor}[defn]{Corollary}

\newcommand{\h}{{\cal H}}

\newcommand{\ltr}{ L^2(\mathbb R) }

\newcommand{\lti}{{\ell}^2(I)}

\newcommand{\si}{S^{-1}}

\newcommand{\mn}{\mathbb N}

\newcommand{\mr}{\mathbb R}

\newcommand{\mz}{\mathbb Z}

\newcommand{\mc}{\mathbb C}

\newcommand{\mts}{ \{E_{mb}T_{na}g \}_{m,n \in \mz}}

\def\bp{{\noindent\bf Proof. \ }}

\def\ep{\hfill$\square$\par\bigskip}

\def\bqs{\begin{equation}}

\def\eqs{\tag*{$\square$}\end{equation}\par\bigskip}

\def\la{\langle}

\def\ra{\rangle}

\def\span{\overline{\text{span}}}

\def\bee{\begin{eqnarray}}

\def\ene{\end{eqnarray}}

\def\bes{\begin{eqnarray*}}

\def\ens{\end{eqnarray*}}

\def\bei{\begin{itemize}}

\def\eni{\end{itemize}}

\def\bt{\begin{thm}}

\def\et{\end{thm}}

\def\bc{\begin{cor}}

\def\ec{\end{cor}}

\def\bpr{\begin{prop}}

\def\epr{\end{prop}}

\def\bl{\begin{lemma}}

\def\el{\end{lemma}}

\def\bd{\begin{defn}}

\def\ed{\end{defn}}

\def\bex{\begin{ex}}

\def\enx{\end{ex}}

\def\bfi{\begin{fig}}

\def\efi{\end{fig}}

\def\eti{\{e_i\}_{i\in I}}

\def\hti{\{h_i\}_{i\in I}}

\def\fti{\{f_i\}_{i\in I}}

\def\gti{\{g_i\}_{i\in I}}

\def\otj{\{\omega_j\}_{j\in I}}

\def\oj{\omega_j}

\def\sui{\sum_{i\in I}}

\def\suj{\sum_{j\in I}}

\newcommand{\seqgri[1]}{\{#1_i\}_{i\in I}}

\newcommand{\seqgrj[1]}{\{#1_j\}_{j\in I}}

\newcommand{\sumgr[1]}{\sum_{#1=1}^\infty}
\def\newin {\,\kern-0.4em\in\kern-0.15em}
\def\newsubset {\kern-0.2em\subset\kern-0.2em}

\def\sumgri{\sum_{i\in I}}
\def\sumgrj{\sum_{j\in I}}

\def\<{\langle}
\def\>{\rangle}

{\nopagebreak\par\noindent\hbox to \hsize{\hrulefill}\par\noindent}

\title{On R-duals and the duality principle in Gabor analysis}

\date{\today}

\author{Diana Stoeva, Ole Christensen}

\begin{document}

\maketitle

\begin{abstract}  The concept of R-duals of a frame was introduced by Casazza, Kutyniok and Lammers in 2004, with the motivation to obtain a general version of the duality principle in Gabor analysis. For tight Gabor frames and Gabor Riesz bases the three authors were actually able to show that the duality principle is a special case of  general results for R-duals.
In this paper we introduce various alternative R-duals, with focus on what we  call R-duals of type II and III.
We show how they are related and provide characterizations of the R-duals of type II and III. In particular, we prove that for tight frames these classes coincide with the R-duals by Casazza et el., which is desirable in the sense that the motivating case of tight Gabor frames already is well covered by these R-duals. On the other hand,  all the introduced types of R-duals generalize the duality principle for larger classes of Gabor frames than just the tight frames and the Riesz bases; in particular, the R-duals of type III cover the duality principle for all Gabor frames.

\end{abstract}

\section{Introduction}
The concept of R-duals of a frame was introduced in the paper \cite{CKL} by Casazza, Kutyniok and Lammers. They showed that the relation between a frame and its R-duals
resembles the known results for the connection between a Gabor system
$\mts$ and  $\{\frac1{\sqrt{ab}} \,
E_{m/a}T_{n/b}g\}_{m,n\in \mz}$ (see the definition below). This lead the three authors
to the natural question whether  $\{\frac1{\sqrt{ab}} \,
E_{m/a}T_{n/b}g\}_{m,n\in \mz}$ can be realized as an R-dual of $\mts;$ while the general question remains open, they were able to confirm this for  Gabor frames $\mts$ that are either tight or Riesz bases.

In this paper we propose various alternative definitions of R-duals, to be called  {\it R-duals of type II, III, IV,} with focus on type II $\&$ III. For this reason we will from now on refer to the R-duals by Casazza $\&$ al. as {\it R-duals of type I.} Among the new types, the R-duals of type II form the smallest set; the R-duals of type III and IV are defined via two steps of relaxation of the conditions. For the case of a tight frame we show that the classes of R-duals of type I, II and III coincide; this is a desirable property because we know that the motivating case of a tight Gabor frames is already well covered by R-duals of type I. Based on a characterization of R-duals of type III we show that for any Gabor frame $\mts$ the sequence
$\{\frac1{\sqrt{ab}} \,
E_{m/a}T_{n/b}g\}_{m,n\in \mz}$ can be realized as an R-dual of type III of $\mts;$ for R-duals of type II  this is at least possible for all integer-oversampled Gabor frames.

In the rest of this introduction we recall the definition of the R-duals of type I, state the main results from \cite{CKL}, and prove a new result that characterizes when a
tight Riesz sequence is an R-dual of a given tight frame with the same frame bound (Proposition \ref{213c}). We also provide the necessary background on Gabor frames, most importantly, the duality principle. Finally, for easy reference we state a few general results about frames and Riesz bases.
In Section \ref{183b} we derive a relationship between $\{\frac1{\sqrt{ab}} \,
E_{m/a}T_{n/b}g\}_{m,n\in \mz}$ and $\mts$ that naturally leads us to the definition of the various types of R-duals.  Section \ref{183c} is devoted to an analysis of the R-duals of type II and the relations between R-duals of type I and II. Section \ref{183d}
presents the properties of R-duals of type III. We provide a characterization of these R-duals (Theorem \ref{213a}), and show that for any given Gabor frame $\mts$ the sequence
$\{\frac1{\sqrt{ab}} \,
E_{m/a}T_{n/b}g\}_{m,n\in \mz}$ can be realized as an R-dual of type III (Corollary \ref{304d}). In the same section we prove that R-duals of type III enjoy most of the properties that make R-duals of type I attractive.

We note that the duality principle and the work on R-duals in \cite{CKL} have trigged quite some interest in various directions. The papers \cite{CKK} by Christensen, Kim, and Kim
and \cite{FS} by Fan and Shen consider R-duals within the original framework of \cite{CKL}, in \cite{FS} using the concept of adjoint systems. In \cite{H} Dutkay, Han and Larson prove that the duality principle extends to any dual pair of projective unitary representations of
countable groups. Finally, X. M. Xiao and Zhu considered an extension of R-duality to Banach spaces in \cite{XZ}, a work that was followed up by  \cite{CXZ} by
Christensen, X. C. Xiao and Zhu.

\subsection{Basic results on frames and Riesz bases}
For easy reference we will collect some of the needed facts about frames and Riesz bases here.  Much more information can be found in the standard monographs, see, e.g.,
\cite{Y,Da2,G2, CBN}.

In the entire paper we let $\h$ denote a separable Hilbert space, with the inner product
$\la \cdot, \cdot \ra$ chosen to be linear in the first entry. We skip the formal definition of a frame and a Riesz basis, which is expected to be well known; we just mention that when we speak about a {\it frame,} it is understood that it is a frame for the entire space $\h.$ In contrast, a {\it frame sequence} is a sequence that is just a frame for the closed span of its elements. We use the same distinction between
a {\it Riesz basis} (which spans $\h$) and a {\it Riesz sequences} (which is a Riesz basis for the closed span of its elements).

It is well known that we can construct a tight frame based on any given frame:

\bl \label{223u} Let $\fti$ be a frame for $\h,$ with frame operator $S.$ Then
$\{S^{-1/2}f_i\}_{i\in I}$ is a tight frame for $\h$ with frame bound 1. If
$\fti$ is a Riesz basis, then $\{S^{-1/2}f_i\}_{i\in I}$ is an orthonormal basis for $\h.$\el

For any given frame there is a natural procedure to construct a Riesz basis with the same frame bounds; see, e.g., \cite{CBN} for a proof of this standard result.

\bl \label{pfr3}
Let $\seqgri[e]$ be any orthonormal basis for $\h$  and $Q:\h \to \h$
a bounded bijective operator. Then the following holds.
\bei
\item[(i)] The sequence $\{Qe_i\}_{i\in I}$ is a Riesz basis with frame
operator $QQ^*$ and optimal bounds $\frac1{||Q^{-1}||^2}, ||Q||^2$.
\item[(ii)] The dual Riesz basis of $\{Qe_i\}_{i\in I}$  is
 $\{(Q^*)^{-1}e_i\}_{i\in I};$ the frame operator for this sequence is $(QQ^*)^{-1}$ and the optimal bounds are $\frac{1}{||Q||^2}, ||Q^{-1}||^2.$
 \eni In particular, if $\fti$ is a frame with frame operator $S$ and optimal bounds $A,B,$ then $\{S^{1/2}e_i\}_{i\in I}$ is a Riesz basis with frame
operator $S$ and optimal bounds $A, B;$ the dual Riesz is $\{S^{-1/2}e_i\}_{i\in I},$ with frame operator $\si$ and optimal bounds $\frac{1}{B}, \frac{1}{A}.$
\el

In case $\fti$ is a frame sequence in $\h$ and
$V:=\span \fti\neq \h,$ the frame operator $S$ can be considered as a bijection on $V.$ The following elementary lemma shows in particular how we can extend $S$ to a bijective and bounded operator on $\h,$ while keeping the norm of the operator and its inverse; we will need this result in the analysis of R-duals of type III.

\bl \label{203a} Let $V$ be a closed subspace of $\h$ and $\Phi:V\to V$ a bounded bijective
operator.
Define an extension of $\Phi$ to
an operator
\bes \widetilde{\Phi}: \h \to \h, \, \widetilde{\Phi}(x_1+x_2):= \Phi x_1 + ||\Phi^{-1}||^{-1} x_2, \, x_1\in V, x_2\in V^\perp.\ens Then $\widetilde{\Phi}$ is bijective and  bounded,  $||\widetilde{\Phi}||=||\Phi||, \,
||\widetilde{\Phi}^{-1}||=||\Phi^{-1}||$, and
\bee \label{303a} \widetilde{\Phi}^{-1}(x_1+x_2)= \Phi^{-1}x_1+ ||\Phi^{-1}||\, x_2, \, \,
x_1\in V, x_2 \in V^\perp\ene If $\Phi$ is self-adjoint, then also $\widetilde{\Phi}$ is self-adjoint.\el

\bp Note that $ ||\Phi^{-1}||^{-1}  \le || \Phi||.$  Thus, writing  $x\in\h$ as $x=x_1 +x_2$ with
 $x_1\in V, x_2\in V^\perp,$
\begin{eqnarray*}
 \| \widetilde{\Phi} x\|^2 &=&\| \Phi x_1\|^2 + \| ||\Phi^{-1}||^{-1}   x_2\|^2\leq || \Phi||^2\, \|x_1\|^2+ ||\Phi^{-1}||^{-2}\|x_2\|^2 \leq ||\Phi||^2\, \|x\|^2. \end{eqnarray*}
From here it is clear that $\Phi$ is bounded and that $||\widetilde{\Phi}||=||\Phi||.$
The other properties are clear by construction. \ep

\subsection{R-duals of type I}
Let $\eti$ and $\hti$ denote orthonormal bases
for $\h,$ and let $\fti$ be any sequence in $\h$ for which $\label{14a} \sui
|\la f_i,e_j\ra|^2 < \infty, \  \forall j\in I.$ In \cite{CKL} the {\it R-dual} of
$\fti$ with respect to the orthonormal bases $\eti$ and $\hti$ is defined as
the sequence $\otj$ given by \bee \label{207a} \oj= \sui \la f_i,
e_j\ra h_i, \ j\in I.\ene As mentioned in the introduction we will from now on refer to
$\otj$ in \eqref{207a} as an {\it R-dual of type I.}
In \cite{CKL} the following connections between the properties of $\fti$ and $\otj$ are proved. Note in particular that an R-dual of a frame is a Riesz sequence:

\bt \label{207b} {\rm \cite{CKL}} Let $\eti$ and $\hti$ denote orthonormal bases
for $\h,$ and let $\fti$ be any sequence in $\h$ for which $\sui
|\la f_i,e_j\ra|^2 < \infty$ for all $j\in I.$ Define the R-dual
$\otj$ of type I  as in \eqref{207a}. Then the following hold:

\bei \item[(i)] For all $i\in I, \, f_i= \suj \la
\oj, h_i\ra e_j,$ i.e., $\fti$ is the R-dual sequence of $\otj$ of type I
w.r.t. the orthonormal bases $\hti$ and $\eti.$
\item[(ii)]  $\fti$ is a frame with bounds $A,B$ if and
only if $\otj$ is a Riesz sequence in $\h$  with bounds $A,B.$
\item[(iii)] Two Bessel sequences $\fti$ and $\gti$ in $\h$ are
dual frames if and only if the associated R-dual sequences $\otj$
and $\{ \gamma_j\}_{j\in I}$ of type I, w.r.t. the same choices of
orthonormal bases $\eti$ and $\hti,$  satisfy \bee
\label{207f} \la \oj, \gamma_k\ra = \delta_{j,k}, \ j,k\in I.\ene
\item[(iv)]  $\seqgrj[\omega]$  is a Riesz basis  if and only if $\seqgri[f]$ is a Riesz basis.
\eni \et
Recall that if $\fti$ is a frame for $\h,$ the {\it preframe operator} or
{\it synthesis operator} is defined by
\bes T: \lti \to \h, \, T\{c_i\}_{i\in I}= \sum_{i\in I} c_if_i.\ens
The following result from \cite{CKL} presents a necessary condition on a sequence $\otj$ to be an R-dual of type I of a given frame $\fti,$ stated in terms of the dimension of the kernel of $T$ and the {\it deficit} of the sequence $\otj:$
\bl \label{213b} {\rm \cite{CKL}} If $\fti$ is a frame with synthesis operator $T$ and
$\otj$ is an R-dual of type I of $\fti,$  then
$ {\rm dim (ker } \, T)= {\rm dim ( span} \otj^\perp.)$\el

The dimension condition in Lemma \ref{213b} will play a crucial role for the various R-duals to be defined in Section \ref{183b}, see Theorem \ref{213a}.
Using Lemma \ref{213b} we can derive a simple characterization of a Riesz sequence $\otj$ being an R-dual of type I of a frame $\fti$ in the tight case:

\bpr \label{213c}
Let $\seqgri[f]$ be a  tight frame for $\h$ and let $\otj$ be a tight Riesz sequence in $\h$ with the same bound. Denote the synthesis operator for $\fti$ by  $T.$ Then  $\otj$ is an R-dual of $\fti$ of type I if and only if
\bee \label{223f} {\rm dim (ker} \, T)= {\rm dim (span}\otj^\perp).\ene
\epr
\bp The necessity of the condition in \eqref{223f} follows from Lemma \ref{213b}. Now
let $\seqgri[f]$ be a  tight frame and $\otj$ be a tight Riesz sequence with the same bound $A,$ and assume that \eqref{223f} holds.
Take an orthonormal basis $\seqgri[e]$ for $\h$
and observe that $\{\frac{1}{\sqrt{A}}\omega_j\}_{j\in I}$ is an orthonormal sequence. Consider the R-dual $\{\nu_j\}_{j\in I}$ of type I of $\fti$ with respect to the orthonormal bases $\seqgri[e]$   and $\seqgri[h]=\seqgri[e]$ , i.e.
$\nu_j := \sumgri \<f_i, e_j\> e_i, \ \, j\in I.$
By Theorem \ref{207b} $\{\nu_j\}_{j\in I}$ is a tight Riesz sequence
with bound $A$ and hence
$\{\frac{1}{\sqrt{A}}\nu_j\}_{j\in I}$ is also an orthonormal sequence.
By Lemma \ref{213b} and \eqref{223f},
\bee \label{223g} {\rm dim (span}\{\nu_j\}_{j\in I}^\perp)={\rm dim (ker} \, T)= {\rm dim(span}\{\omega_j\}_{j\in I}^\perp).\ene

In case ${\rm span}\{\nu_j\}_{j\in I}^\perp ={\rm span}\{\omega_j\}_{j\in I}^\perp =\{0\},$ the orthonormality of the sequences $\{\frac{1}{\sqrt{A}}\nu_j\}_{j\in I},
\{\frac{1}{\sqrt{A}}\omega_j\}_{j\in I},$ implies that we can define a unitary  operator $U: \h \to \h$ by
$
U\nu_j :=\omega_j, \ j\in I;
$ and in case ${\rm span}\{\nu_j\}_{j\in I}^\perp \neq \{0\},$ letting
$\{\phi_j\}_{j\in J}$ and $\{\psi_j\}_{j\in J}$ be orthonormal bases for ${\rm span}\{\nu_j\}_{j\in I}^\perp$ and ${\rm span}\{\omega_j\}_{j\in I}^\perp, $
respectively, \eqref{223g} has the consequence that we can define a unitary operator $U: \h\to \h$ by
$U\nu_j :=\omega_j, j\in I,
\mbox{ and }
U\phi_j := \psi_j, \ j\in J.
$
In both cases,
\begin{equation*}
 \omega_j = U\nu_j
= U\sumgri \<f_i, e_j\> e_i
=
\sumgri \<f_i, e_j\>  U e_i, \ j\in I,
\end{equation*}
which shows that $\otj$ is an R-dual of $\fti$ of type I. \ep

Note that in the non-tight case
\eqref{223f} does not imply that $\otj$ is an R-dual of $\fti$ of type I:

\bex \label{293g}
Let $\{z_i\}_{i=1}^\infty$ be an orthonormal basis for $\h$. Consider the frame (actually a  Riesz basis)
$ \{f_i\}_{i=1}^\infty=\{\sqrt{2}z_1, z_2, \sqrt{2}z_3, \sqrt{2} z_4, \sqrt{2} z_5, \sqrt{2} z_6, \ldots\}$
and the Riesz basis
$\{g_j\}_{j=1}^\infty=\{\sqrt{2}z_1, z_2, z_3, z_4, z_5, \ldots\}.$
The sequences $ \{f_i\}_{i=1}^\infty$ and $\{g_j\}_{j=1}^\infty$ both have the
optimal  bounds $A=1$, $B=2$ and \eqref{223f} holds, but $ \{g_j\}_{i=j}^\infty$ is not an R-dual of $\{f_i\}_{i=1}^\infty$ of type I. In fact,
assume that there exist orthonormal bases $\{e_i\}_{i=1}^\infty$ and $\{h_i\}_{i=1}^\infty$
for $\h$ so that
$g_j=\sumgr[i] \<f_i, e_j\> h_i, \ \forall j\in\mn.$
Then
$
z_{j} = g_j=\sumgr[i] \<f_i, e_j\> h_i, \ \mbox{for every }  j\geq 2, j\in\mn,
$
which implies that $\< z_{j} , h_i\> =\<f_i, e_j\>$ for all $i,j\in\mn, j\geq 2$.
Thus
\begin{eqnarray*}
 1 &=& \|z_{j} \|^2 = \sum_{i=1}^\infty |\<z_{j} , h_i\>|^2
 = \sum_{i=1}^\infty |\<f_i, e_j\>|^2=
\\
 &=&  \sum_{i=1}^\infty |\<z_i, e_j\>|^2 + \sum_{i\in\{1,3,4,5,\ldots\}} |\<z_i, e_j\>|^2  =
 1 + \sum_{i\in\{1,3,4,5,\ldots\}}  |\<z_i, e_j\>|^2.
 \end{eqnarray*}
It follows that $z_1 \perp e_j$ and $z_3\perp e_j$ for all $j\in \mn, j\geq 2$,
which is a contradiction. Thus $ \{g_j\}_{j=1}^\infty$ is not an R-dual of $\{f_i\}_{i=1}^\infty$ of type I.\ep \enx

\subsection{Gabor analysis}
For parameters $a,b\in \mr,$ define the operators $T_a$ and $E_b$ on $\ltr$ by
$T_af(x)=f(x-a)$ and $E_bf(x)=e^{2\pi ibx}f(x),$ respectively. The {\it Gabor system}
generated by a fixed function $g\in \ltr$ and some $a,b>0$ is the collection of functions $\mts.$ Remember that when we speak about $\mts$ being a frame, it is understood that we mean a frame for the space $\ltr.$ For a detailed discussion of the role of systems
$\mts$ in time-frequency analysis we refer to the monograph \cite{G2}; their frame properties are in focus in the monograph \cite{CBN}.

One of the deepest results in Gabor analysis is the {\it duality principle}, which  was discovered at the same time by Janssen \cite{Jan5}, Daubechies,
Landau, and Landau \cite{DLL}, and Ron and Shen \cite{RoSh1}:

\begin{thm} \label{1505d1} {\rm \cite{DLL, Jan5, RoSh1}} Let $g\in \ltr$ and $a,b>0$ be given.
Then the Gabor system $\mts$ is a frame with bounds $A,B$ if and only if $\{\frac1{\sqrt{ab}} \,
E_{m/a}T_{n/b}g\}_{m,n\in \mz}$ is a Riesz sequence with bounds $A,B.$
\end{thm}

The similarity between Theorem \ref{207b}(ii) and Theorem \ref{1505d1}  leads to the obvious question whether there is a connection between the duality principle and the R-duals of type I. In \cite{CKL} it was shown that at least for two important special cases of Gabor frames
$\mts,$ the sequence $\{\frac1{\sqrt{ab}} \,
E_{m/a}T_{n/b}g\}_{m,n\in \mz}$ is actually an R-dual of $\mts:$

\bt \label{22a} {\rm \cite{CKL}} Assume that the frame $\mts$ is either tight  or a Riesz basis.
Then $\{\frac1{\sqrt{ab}} \,
E_{m/a}T_{n/b}g\}_{m,n\in \mz}$ can be realized as an R-dual of type I of
$\mts.$\et

It is still an open question whether Theorem \ref{22a} holds for other classes of Gabor frames. We will not provide any new results about this; instead, our purpose is to introduce other types of R-duals that cover the duality principle for larger classes of Gabor frames than the two classes in Theorem \ref{22a}. The R-duals of type III (to be defined in Section \ref{183b} and analyzed in Section
\ref{183d}) turns out to cover all Gabor frames, and also keep essential parts of the properties  R-duals of type I, as stated in Theorem \ref{207b}.

Let us end this section with a few results from Gabor analysis that will be used repeatedly. The first one is due to Balan, Casazza, and Heil. Recall that a frame $\fti$ has {\it infinite excess} if infinitely many elements can be removed while the remaining sequence is still a frame.
Also, the {\it deficit} of a sequence $\fti$ in $\h$ is the number
${\rm dim}( \span \fti^\perp).$

\bl {\rm \cite{BCH}} \label{293a} Let $g\in \ltr$ and $a,b>0$ be given. Then the following hold.
\bei
\item[(i)] If $ab<1$ and $\mts$ is a frame, then $\mts$ has infinite excess.
\item[(ii)] If $ab>1,$ then $\mts$ has infinite deficit.\eni\el

The following well-known result will be used at several places. Recall that for any frame, the frame operator $S$ is a positive bounded operator, and thus has a unique positive square-root $S^{1/2}.$

\bl \label{223p} Let $\mts$ be a Bessel sequence in $\ltr,$ with frame operator $S.$
Then $S$ commutes with the operators $E_{mb}T_{na}, \, m,n\in \mz.$ If $\mts$ is a frame, then also the operators $S^{1/2}$ and $S^{-1/2}= (S^{-1})^{1/2}$
commute with the operators $E_{mb}T_{na}, \, m,n\in \mz.$ \el
Recall that a Gabor frame $\mts$ is said to be {\it integer-oversampled} if  $ab=1/K$ for some $K\in \mn.$ For our purpose the importance of integer-oversampled Gabor frames lies in the fact that for such frames the operators $E_{m/b}T_{n/a}, m,n\in \mz,$ form a {\it subclass} of $E_{mb}T_{na}, \, m,n\in \mz$ and therefore commute with the frame operator
of $\mts$  and its various powers.

\section{Towards generalized R-duals} \label{183b}
We have already noticed that for Gabor frames $\mts$ that are either tight or Riesz bases, the sequence $\{\frac1{\sqrt{ab}} \,
E_{m/a}T_{n/b}g\}_{m,n\in \mz}$ can be realized as an R-dual of
$\mts$ of type I. We will now show that we can cover a larger class of Gabor frames by replacing
the orthonormal bases in the definition of the R-dual by other types of sequences.

\bpr \label{103b} Assume that $\mts$ is either a tight Gabor frame or an integer-oversampled Gabor frame, with frame operator $S.$ Then there exist Riesz bases $\{x_{m,n}\}_{m,n\in \mz}, \{y_{m,n}\}_{m,n\in \mz}$ for $\ltr$ such that
\bee \label{103a} \frac1{\sqrt{ab}} \,
E_{m/a}T_{n/b}g = \sum_{m^\prime,n^\prime\in \mz} \la E_{m^\prime b}T_{n^\prime a}g, x_{m,n}\ra y_{m^\prime, n^\prime}.\ene More explicitly, there exist orthonormal bases $\{e_{m,n}\}_{m,n\in \mz}$ and $\{h_{m,n}\}_{m,n\in \mz}$ such that
\eqref{103a} holds with $x_{m,n}= S^{-1/2}e_{m,n}, \, y_{m,n}= S^{1/2}h_{m,n}.$\epr

\bp Let us first consider the case where $\mts$ is an integer-oversampled Gabor frame. Then
$\{E_{mb}T_{na} S^{-1/2}g\}_{m,n\in \mz}$ is a tight frame for $\ltr$ with frame bound 1. Let $\{e_{m,n}\}_{m,n\in \mz}$ be an orthonormal basis for $\ltr$ consisting of bounded and compactly supported functions. Then, according to the proof of Theorem \ref{22a}(ii) in \cite{CKL} there exists a unitary operator $U: \ltr\to \ltr$ such that
\bee \label{193a} \frac1{\sqrt{ab}} \,
E_{m/a}T_{n/b} S^{-1/2} g = \sum_{m^\prime,n^\prime\in \mz} \la E_{m^\prime b}T_{n^\prime a}S^{-1/2}g, e_{m,n}\ra Ue_{m^\prime, n^\prime}.\ene
Since $\mts$ is a frame, the operator $S^{-1/2}$ commutes with the operators
$E_{mb}T_{na}$ for all $m,n\in \mz.$ Since $ab=1/K$ for some $K\in \mn,$ the operators $
E_{m/a}T_{n/b}$ form a subclass of the operators $E_{mb}T_{na}, \, m,n\in \mz,$ and therefore \eqref{193a}
implies that
\bes \frac1{\sqrt{ab}} \,
S^{-1/2}E_{m/a}T_{n/b}  g = \sum_{m^\prime,n^\prime\in \mz} \la S^{-1/2}E_{m^\prime b}T_{n^\prime a}g, e_{m,n}\ra Ue_{m^\prime, n^\prime}.\ens Thus
$\frac1{\sqrt{ab}} \,
E_{m/a}T_{n/b}  g = \sum_{m^\prime,n^\prime\in \mz} \la E_{m^\prime b}T_{n^\prime a}g, S^{-1/2}e_{m,n}\ra S^{1/2}Ue_{m^\prime, n^\prime},$ which proves \eqref{103a},
with sequences $x_{m,n}$ and $y_{m,n}$ as stated in the proposition.

Let us now assume that $\mts$ is a tight frame, with frame bound $A.$ Then
$S=AI.$ By Theorem \ref{22a} there exist orthonormal bases  $\{e_{m,n}\}_{m,n\in \mz}$
and
$\{h_{m,n}\}_{m,n\in \mz}$ such that
\bee \label{404a} \frac1{\sqrt{ab}} \,
E_{m/a}T_{n/b}g = \sum_{m^\prime,n^\prime\in \mz} \la E_{m^\prime b}T_{n^\prime a}g, e_{m,n}\ra h_{m^\prime, n^\prime}.\ene Using that $S^{1/2}= A^{1/2}I$ and
$S^{-1/2}= A^{-1/2}I,$ \eqref{404a} can clearly be written on the form \eqref{103a}, as claimed.
\ep

Proposition \ref{103b} leads to several natural ways of alternative definitions of R-duals. We will
call the new types for  {\it R-duals of type II, III, IV,} respectively. Note that at the moment we have not motivated the definition of R-duals of type III; the reason for the definition will become clear in Section \ref{183d}, where we prove a statement of similar spirit as Proposition \ref{103b} but without the assumption on $\mts$ being integer-oversampled. The resulting expansion places the general Gabor case within the framework of R-duals of type III.

\begin{defn} \label{253a} Let $\seqgri[f]$ be a frame for $\h$ with frame operator $S.$
\bei
\item[(i)] Let $\seqgri[e]$ and $\seqgri[h]$ denote orthonormal bases for $\h$. The R-dual of type II of $\seqgri[f]$ with respect to $\seqgri[e]$ and $\seqgri[h]$ is the sequence  $\seqgrj[\omega]$ given by
\begin{equation}\label{grdual}
 \omega_j=\sumgri \<  f_i, S^{-1/2} e_j\> S^{1/2} h_i, \ j\in I.
 \end{equation}
 \item[(ii)]
Let $\seqgri[e]$ and $\seqgri[h]$ denote orthonormal bases for $\h$
and $Q:\h\to\h$ be a bounded bijective
operator with  $\|Q\|\leq \sqrt{||S||}$ and  $\|Q^{-1}\|\leq \sqrt{|| \si ||}$.
The R-dual of type III of $\seqgri[f]$ with respect to the triplet ($\seqgri[e]$, $\seqgri[h]$, $Q$), is the sequence  $\seqgrj[\omega]$ defined by
\begin{equation}\label{rfdualGab}
\omega_j:=\sumgr[i] \< S^{-1/2}f_i, e_j\> Qh_i.
\end{equation}
 \item[(iii)] Let $\seqgri[e]$ and $\seqgri[h]$ denote Riesz bases for $\h$. The R-dual of type IV of $\seqgri[f]$ with respect to $\seqgri[e]$ and $\seqgri[h]$ is the sequence  $\seqgrj[\omega]$ given by
\begin{equation}\label{grdualp}
 \omega_j=\sumgri \<  f_i, e_j\> h_i, \ j\in I.
 \end{equation}
\eni
\end{defn}

Note that in the definition of R-duals of type III the action of the operator
$S^{-1/2}$ is written on $f_i$ instead of $e_j.$ This implies that R-duals of type III
also are defined if we only assume that $\fti$ is a frame sequence and interpret the frame operator as a bijection on
$\span \fti.$ In contrast, R-duals of type II only exist when $\fti$ is a frame for $\h.$

The conditions on the operator $Q$ in Definition \ref{253a} (ii) means that $\{Qh_i\}_{i\in I}$ is a Riesz basis for $\h$ with bounds $||\si||^{-1}, ||S||,$ which according to Lemma \ref{pfr3} are the optimal bounds for the sequence $\{S^{1/2} h_i\}_{i\in I}.$ Thus, the R-duals of type II are contained in the class of R-duals of type III. It is obvious that the R-duals of type III are contained in the class of R-duals of type IV. It is also clear that R-duals of type I are contained in the class of R-duals of type IV. Example \ref{exb1}, Example \ref{293b},
and Theorem \ref{213a} provide further information about the relationship between the various classes; in order for the reader to get a quick overview of what to come, we summarize these results in Figure 1.
We also note already now that Proposition \ref{304a} and Proposition \ref{13coincide}
will show that for a tight frame, the classes of R-duals of type I, II and III coincide.

\begin{figure}[hbt]\label{fig1}
\begin{center}
\mbox{
{\includegraphics
[width=60mm,height=32mm]
{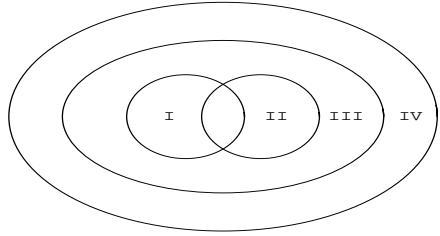}}}
\hspace{2.2in}\caption{The relationship between R-duals of type I, II, III, and IV.}
\end{center}
\end{figure}

A general analogue of R-duals of type IV have been considered in Banach spaces in the papers \cite{XZ} and \cite{CXZ}.  In particular, the results by Xiao and Zhu in \cite{XZ} imply that in the framework of Definition \ref{253a} (iii),
$f_i= \sum_{j\in I} \la \omega_j, \widetilde{h_i}\ra \widetilde{e_j},$ where $\{\widetilde{e_i}\}_{i\in I}$ and $\{\widetilde{h_i}\}_{i\in I}$
are the dual Riesz bases of $\eti$ and $\hti,$ respectively; furthermore, if $\fti$ is a frame then $\otj$ is a Riesz sequence, with  a lower bound that is the multiple of the lower bound of $\fti$ and the lower bounds for the two Riesz sequences $\eti, \hti$  (and a similar result for the upper bound).  These properties are getting too far away from the properties we know for R-duals of type I, so we will not discuss type IV in more detail in this paper.

\section{R-duals of type II} \label{183c}

The purpose of this section is to state some of the most important properties of R-duals of type II and to compare the properties of the R-duals of type I and II. Other results
are on Section \ref{183d}, which cover the more general case of R-duals of type III. Also note that Theorem \ref{213a} contains a characterization of the R-duals of type II.

First we notice that for tight frames, the classes of R-duals of type I and II coincide:

\bpr \label{304a} Assume that $\fti$ is a tight frame for $\h.$ Then the classes of R-duals of type I and II coincide. \epr
\bp Since the frame operator $S=AI$ for some $A>0,$ $S^{1/2}= A^{1/2}I$ and
$S^{-1/2}= A^{-1/2}I.$ Thus, for any orthonormal bases $\eti$ and $\hti,$
we have $\sum_{i\in I} \la f_i, e_j\ra h_i =   \sum_{i\in I} \la f_i, S^{-1/2}e_j\ra
S^{1/2}h_i,$ which proves the result. \ep

In general the R-duals of type I and II constitute two different classes. The following example exhibits an R-dual of type II of a specific frame, which is not an R-dual of type I; in Example \ref{293b} we will use the same construction to find an R-dual of type I which is not an R-dual of type II.

\begin{ex} \label{exb1}
Let $\{e_i\}_{i=1}^\infty$ be an orthonormal basis for $\h$.
Consider the frame $\{f_i\}_{i=1}^\infty=\{e_1, e_1, e_2, e_3, e_4, \ldots\},$
which has optimal bounds $1,2.$  Note that $S e_1=2e_1, Se_i=e_i, e\geq 2$.
A simple calculation of the R-dual of type II of $\{f_i\}_{i=1}^\infty$ with respect to $\{e_i\}_{i=1}^\infty$ and $\{h_i\}_{i=1}^\infty=\{e_i\}_{i=1}^\infty$ shows that
\\ $ \{\omega_j\}_{j=1}^\infty= \{ e_1+ \frac1{\sqrt{2}}\, e_2, e_3, e_4, \dots\}.$
In particular, $\{\omega_j\}_{j=1}^\infty$ has the optimal bounds $1, 3/2.$  By Theorem \ref{207b} the R-duals of type I have the same optimal bounds as the given frame, which shows that $\{\omega_j\}_{j=1}^\infty$ can not be an R-dual of type I of $\{f_i\}_{i=1}^\infty.$\ep \enx

It follows immediately from Proposition \ref{103b} that for any  Gabor frame $\mts$ that is either tight or integer-oversampled, the sequence  $\{\frac1{\sqrt{ab}} \,
E_{m/a}T_{n/b}g\}_{m,n\in \mz}$ can be realized as an R-dual of type II. In the general setting of frames in Hilbert spaces we now show that R-duals of type II enjoy many of the properties that make R-duals of type I attractive, see Theorem \ref{207b}. This is the content of the following proposition, as well as Proposition \ref{304b} and Theorem \ref{213a}
in the next section.

\begin{prop} \label{generalrdual}
\label{rrdual}
Let $\seqgri[f]$ be a frame for $\h$ with frame operator $S$ and let $\seqgrj[\omega]$ be the R-dual of type II  of $\seqgri[f]$  with respect to some orthonormal bases $\seqgri[e]$ and $\seqgri[h]$.
The following statements hold.
\begin{itemize} \item[(i)]
$ f_i =  \sumgrj \<\omega_j, S^{-1/2}h_i\> S^{1/2}e_j$, $\forall i\in I$.

\item[(ii)]   Let $\{\gamma_j\}_{j\in I}$ denote the R-dual of type II of the
canonical dual frame $\{\si f_i \}_{i\in I}$ with respect to
$\seqgri[e]$ and $\seqgri[h].$ Then $\otj$ and $\{\gamma_j\}_{j\in I}$ are biorthogonal.

\end{itemize}
\end{prop}

\bp (i) By Lemma \ref{pfr3} and \eqref{grdual},  $\<\omega_j, S^{-1/2}h_i\>  = \<f_i, S^{-1/2}e_j\>$ for every $i, j\in I;$ thus, expanding $f_i$ with respect to the dual pair of Riesz bases $\{ S^{-1/2}e_j\}_{j\in I}, \{ S^{1/2}e_j\}_{j\in I}$ yields that
$$ f_i = \sumgrj \<f_i, S^{-1/2}e_j\> S^{1/2}e_j = \sumgrj \<\omega_j, S^{-1/2}h_i\> S^{1/2}e_j, \ \forall i\in I.$$

\noindent(ii) The frame operator of $\{\si f_i \}_{i\in I}$ is $\si,$ so by definition \\
$\gamma_n= \sum_{k\in I} \<  \si f_k, S^{1/2} e_n\> S^{-1/2} h_k, \, n\in I.$
Therefore, for $j,n\in I,$
\begin{eqnarray*}
\<\omega_j,\gamma_n\>
&=& \<\sumgri \<  f_i, S^{-1/2} e_j\> S^{1/2} h_i, \sum_{k\in I} \<  \si f_k, S^{1/2} e_n\> S^{-1/2} h_k\>\\
&=& \sumgri  \sum_{k\in I}  \<  f_i, S^{-1/2} e_j\> \<  S^{1/2} e_n, \si f_k\>  \<S^{1/2} h_i,S^{-1/2} h_k\>.\ens
By Lemma \ref{pfr3} $\{S^{1/2}h_i\}_{i\in I}$ and
$\{S^{-1/2}h_i\}_{i\in I}$ are biorthogonal; thus,
\bes \<\omega_j,\gamma_n\>
&=& \sumgri  \<  f_i, S^{-1/2} e_j\> \<  S^{1/2} e_n, \si f_i\>
=  \<  \sumgri  \<  S^{1/2} e_n, \si f_i\>  f_i, S^{-1/2} e_j\> \\
&= & \<S^{1/2} e_n, S^{-1/2} e_j\> = \delta_{n,j},
\end{eqnarray*} which completes the proof.
\ep

\section{R-duals of type III} \label{183d}
The main result in this section is a characterization of the R-duals of type II and III
(Theorem \ref{213a}). As a consequence we show that for any Gabor frame $\mts,$ the
sequence $\{\frac{1}{\sqrt{ab}} E_{m/a}T_{n/b} g\}_{m,n\in \mz}$ can be
realized as an R-dual of type III (Corollary \ref{304d}). We will also derive some key properties of R-duals of type III and relate them to the properties of R-duals of type I.

We first note that R-duals of type III have an obvious characterization in terms of R-duals of type I. We leave the proof to the reader:
\bl \label{304g} Let $\fti$ be a frame sequence with frame operator $S,$ $\otj$ a Riesz sequence, and $Q: \h \to \h$ a bounded bijective operator with $\|Q\|\leq \sqrt{||S||}$ and  $\|Q^{-1}\|\leq \sqrt{|| \si ||}$. Then the following are equivalent:
\bei \item[(i)] $\otj$ is an R-dual of type III of $\fti$  w.r.t. the operator $Q;$
 \item[(ii)] $\{Q^{-1}\omega_j\}_{j\in I}$ is an R-dual  of type I of $\{S^{-1/2}f_i\}_{i\in I} $.
\eni    In particular, if $\otj$ is an R-dual of $\fti$ of type III w.r.t. $Q,$ then  $\{Q^{-1}\omega_j\}_{j\in I}$
is an orthonormal sequence. \el

We have already seen that for a tight frame, the classes of R-duals of type I and II coincide. In this special case the type III duals give the same class:

\bpr \label{13coincide}
Assume that $\fti$ is a tight frame for $\h.$ Then the classes of R-duals of type I and III coincide. \epr
\bp  Denote the frame operator  for $\fti$ by $S$. Then $S=AI$ for some $A>0$.
Let  $\eti$ and $\hti$ be any orthonormal bases for $\h$.

Take a bounded bijective operator $Q$ which has the property $\|Q\|\leq \sqrt{||S||}= A^{1/2}$ and  $\|Q^{-1}\|\leq \sqrt{|| \si ||}= A^{-1/2}$. By Lemma \ref{pfr3}
$\{Qh_i\}_{i\in I}$ is a tight Riesz basis with   bound $A,$ which implies that it has the form $\{A^{1/2} u_i\}_{i\in I}$ for some orthonormal basis $\{u_i\}_{i\in I}$. Thus the R-dual of $\fti$ of type III  with respect to the triplet ($\seqgri[e]$, $\seqgri[h]$, $Q$) is
 \bes \omega_j= \sum_{i\in I} \la S^{-1/2} f_i, e_j\ra Q h_i =
\sum_{i\in I} \la A^{-1/2} f_i, e_j\ra A^{1/2} u_i
=   \sum_{i\in I} \la f_i, e_j\ra u_i, \ j\in I,\ens
which is an R-dual  of type I of $\fti$.

Now consider the R-dual of $\fti$ of type I  w.r.t. $\seqgri[e]$ and $\seqgri[h],$
\bes \nu_j = \sum_{i\in I} \la f_i, e_j\ra h_i =  \sum_{i\in I} \la A^{-1/2} f_i, e_j\ra A^{1/2} h_i
=\sum_{i\in I} \la S^{-1/2} f_i, e_j\ra S^{1/2} h_i, \ j\in I; \ens
this is clearly an R-dual  of type III of $\fti$.
\ep

For R-duals of type III we now state the analogue of
the properties in Theorem \ref{207b} (ii) and (iv).

\bpr \label{304b} Let $\seqgri[f]$ be a frame  sequence and  $\seqgri[\omega]$ an R-dual of $\fti$ of type III. Then  the following hold.

\bei\item[(i)] $\seqgri[f]$ is a frame if and only if $\otj$ is a Riesz sequence; in the affirmative case the bounds for $\fti$ are also bounds for $\otj.$
\item[(ii)] $\seqgri[f]$ is a Riesz sequence if and only if $\otj$ is a frame; in the affirmative case the bounds
    for $\fti$ are also bounds for $\otj.$
\item[(iii)] $\seqgrj[\omega]$  is a Riesz basis if and only if $\seqgri[f]$ is a Riesz basis. \eni \epr

\bp (i) Assume first that $\fti$ is a frame. Lemma \ref{pfr3} and \eqref{rfdualGab} yield that for any finite scalar sequence $\{c_j\}$,
\begin{eqnarray*}
\| \sum_j c_j \omega_j \|&=& \| \sum_j c_j \sumgri \<f_i, S^{-1/2} e_j\> Q h_i\|
= \| \sumgri  \<S^{-1/2} f_i, \sum_j \overline{c_j} e_j\>  Q h_i\| \\
 &\leq& ||Q|| \, \|\{\<S^{-1/2}f_i, \sum_j \overline{c_j} e_j\>\}_{i\in I} \|_{\ell^2}
= \sqrt{||S||} \, \| \sum_j \overline{c_j} e_j \|_{\h}\\ &=& \sqrt{||S||} \, \| \{c_j\} \|_{\ell^2}.
\end{eqnarray*}
The lower bound is proved in the same way.

On the other hand assume that $\otj$ is a Riesz sequence. Then the sequence
$\{\nu_j\}_{j\in I}$ given by
$\nu_j= Q^{-1}\omega_j= \sum_{i\in I} \la S^{-1/2} f_i,e_j\ra h_i$ is also a
Riesz sequence, which, by Theorem \ref{207b} implies that $\{S^{-1/2}f_i\}_{i\in I}$ is a frame, in particular, that $\{S^{-1/2}f_i\}_{i\in I}$ is total in $\h.$ Since
the frame operator and its powers are bijections on $\overline{span} \fti,$ this implies that $\overline{span} \fti=\h,$ i.e., $\fti$ is a frame for $\h.$

For the proof of (ii), assume that  $\fti$ is a Riesz sequence. By Lemma \ref{304g},
$\{Q^{-1}\omega_j\}_{j\in I}$ is an R-dual of type I of the 1-tight Riesz sequence
$\{S^{-1/2}f_i\}_{i\in I}$, which by Theorem \ref{207b} implies that $\{Q^{-1}\omega_j\}_{j\in I}$ is a tight frame for $\h,$ with bound 1. It follows that $\seqgrj[\omega]$ is  a frame for $\h$ with optimal bounds $ 1/\|Q^{-1}\|^2\geq 1/\|S^{-1}\|$, $\|Q\|^2\leq \|S\|$.

Now assume that $\seqgrj[\omega]$ is  a frame for $\h$. Then the frame $\{Q^{-1}\omega_j\}_{j\in I}$ is an R-dual of type I of
$\{S^{-1/2}f_i\}_{i\in I}$, which by Theorem \ref{207b} implies that $\{S^{-1/2}f_i\}_{i\in I}$  is a Riesz sequence. Therefore, $\fti$ is a Riesz sequence.

In order to prove (iii), assume that $\seqgri[f]$ is a Riesz basis for $\h$. By (ii), $\seqgrj[\omega]$ is a Riesz sequence. It remains to prove the completeness of $\seqgrj[\omega]$.
Let $x\in\h$ and $\<x,\omega_j\>=0$ for every  $j\in I$. Then for   every  $j\in I$,
$$0=\sumgri \<S^{-1/2} e_j,f_i\> \<x, Q h_i\> = \<e_j, \sumgri  \<Q h_i,x\> S^{-1/2} f_i\>,$$
which implies that $\sumgri  \<Q h_i,x\> S^{-1/2} f_i=0$. Therefore, $\<Q h_i, x\>=0$ for every $i\in I$, which implies that $x=0$. Thus $\otj$ is a Riesz basis.

Now assume that $\seqgrj[\omega]$  is a Riesz basis for $\h$.
Consider the sequence $\seqgri[\nu]$ given by $\nu_i= \sumgrj \<\widetilde{\omega}_j, Qh_i\> S^{-1/2}e_j$,  $i\in I$, where $\{\widetilde{\omega}_j\}_{j\in I}$ is the canonical dual of $\seqgrj[\omega].$ We will now show that
$\seqgrj[\nu]$ is biorthogonal to $\seqgri[f].$ Note that by (i) we know that
$\fti$ is a frame. We will now use a representation of $f_i,$ to be proved in Theorem \ref{273b} (i).
Using that \\ $f_i =  \sumgrj \<\omega_j, (Q^*)^{-1} h_i\> S^{1/2}e_j ,$   for every $i,k\in I$ we have
\begin{eqnarray*}
\<\nu_i, f_k\>
 & =&
\sumgrj \<\widetilde{\omega}_j, Qh_i\> \< S^{-1/2}e_j, \sum_{\ell\in I} \<\omega_\ell, (Q^*)^{-1}h_k\> S^{1/2}e_\ell\> \\
& =&
\sumgrj \<\widetilde{\omega}_j, Qh_i\> \sum_{\ell\in I} \< (Q^*)^{-1}h_k, \omega_\ell\> \< S^{-1/2}e_j,  S^{1/2}e_\ell\>.
\end{eqnarray*}
Using Lemma \ref{pfr3}, it now follows that
\begin{eqnarray*}
\<\nu_i, f_k\>
 & =&
\sumgrj \<\widetilde{\omega}_j, Qh_i\> \< (Q^*)^{-1}h_k, \omega_j\>
= \< \sumgrj \< (Q^*)^{-1}h_k, \omega_j\> \widetilde{\omega}_j, Qh_i\> \\
 & =& \< (Q^*)^{-1}h_k,  Qh_i\> =\delta_{k,i}.
\end{eqnarray*} The fact that $\fti$ is a frame and has a  biorthogonal sequence now implies that
$\seqgri[f]$ is  a Riesz basis. \ep

Note that in contrast with the statement for R-duals of type I (Theorem \ref{207b}), Proposition \ref{304b} (i)
only states that frame bounds for $\fti$ are also bounds for
$\otj.$ This statement can not be strengthened, as demonstrated by Example \ref{exb1}, even if $\otj$ is an R-dual of type III with respect to $Q=S^{1/2}$ which satisfies $\|Q\|=\|S\|^{1/2}$ and   $\|Q^{-1}\|= \sqrt{|| \si ||}$. However, if
$\otj$ is assumed to be an R-dual of type III of $\fti$ with respect to $Q$ satisfying $\|Q\|=\|S\|^{1/2}$ and $\|Q^{-1}\|= \sqrt{|| \si ||}$,
then in Proposition \ref{304b} (ii) and (iii), $\otj$ keeps the optimal frame bounds of $\fti$.

We now prove a characterization of R-duals of type II and III:

\bt \label{213a}
Let $\seqgri[f]$ be a  frame for $\h,$ let $\otj$ be a Riesz sequence in $\h$ and assume that the bounds of $\fti$ are also bounds for $\otj.$ Denote the synthesis operator for $\fti$ by  $T.$ Then the following hold.
\bei
\item[(i)] $\otj$ is an R-dual  of type II of $\fti$   if and only if $\{ S^{-1/2} \omega_j\}_{j\in I}$ is an orthonormal system and
\bee \label{273a} {\rm dim (ker} \, T)= {\rm dim (span}\{\omega_j\}_{j\in I}^\perp) \ene
\item[(ii)]  $\otj$ is an R-dual  of type III of $\fti$   if and only if \eqref{273a} holds.
    \item[(iii)] The class of type I duals of $\fti$ is contained in the class of type III duals.\eni
\et
\bp We first prove (ii). Assume that $\otj$ is an R-dual of type III of $\fti$ with respect to some orthonormal bases and some bounded bijective  operator $Q$.
By Lemma \ref{304g} this implies that $\{Q^{-1}  \omega_j\}_{j\in I}$ is an R-dual of type I of $\{S^{-1/2} f_i\}_{i\in I}.$ Since the synthesis operator for $\{S^{-1/2} f_i\}_{i\in I}$ equals $S^{-1/2}T,$
its kernel equals the kernel of $T;$ thus,
Lemma \ref{213b} implies that
\bes {\rm dim (ker} \, T)={\rm dim (span}\{Q^{-1}\omega_j\}_{j\in I}^\perp) ={\rm dim (span}\{\omega_j\}_{j\in I}^\perp), \ens i.e., \eqref{273a} holds. On the other hand, assume that \eqref{273a} holds, and denote the frame operator for $\otj$ by $S_{\Omega}.$ Then $\{ S_{\Omega}^{-1/2}\omega_j\}_{j\in I}$ is a tight Riesz sequence, and
$\{S^{-1/2} f_i\}_{i\in I}$ is a tight frame, both of them with bound 1. Considering
again the
synthesis operator $S^{-1/2}T $ for $\{S^{-1/2} f_i\}_{i\in I},$ the assumption \eqref{273a} implies that
\bes {\rm dim (ker} \, (S^{-1/2}T))={\rm dim (ker} \, T)={\rm dim (span}\{\omega_j\}_{j\in I}^\perp)= {\rm dim (span}\{S_{\Omega}^{-1/2}\omega_j\}_{j\in I}^\perp) .\ens Thus, by Proposition \ref{213c} $\{S_{\Omega}^{-1/2} \omega_j\}_{j\in I}$ is an R-dual of
type I of $\{S^{-1/2} f_i\}_{i\in I},$ i.e., there exist orthonormal bases $\seqgri[e]$ and $\seqgri[h]$   such that
\begin{equation} \label{303g}
S_{\Omega}^{-1/2}  \omega_j= \sumgri \< f_i, S^{-1/2}e_j\>  h_i, \ j\in I.
\end{equation}
Now consider the extension $\widetilde{S_\Omega^{1/2}}$ of $S_\Omega^{1/2}$ to an operator on $\h,$ as in Lemma \ref{203a}.
Then
$||\widetilde{S_\Omega^{/2}} ||= ||S_{\Omega}^{1/2}|| \leq||S||^{1/2}
\ \mbox{and} \
 || (\widetilde{S_\Omega^{/2}})^{-1}||= ||S_{\Omega}^{-1/2}|| \leq ||S^{-1}||^{1/2}.$
Applying the operator $\widetilde{{S}_\Omega^{1/2}}$ to the  representation \eqref{303g} we get
$$   \omega_j  = \sumgri \<f_i, S^{-1/2}e_j\>\widetilde{{S}_\Omega^{1/2}} h_i,
$$ which shows that $\otj$ is an R-dual  of type III of $\fti.$

For the proof of (i), assume that $\otj$ is an R-dual of $\fti$ of type II.
Since R-duals of type II are special cases of the R-duals of type III, the above argument shows that \eqref{273a} holds. Also, by definition
 $\{S^{-1^/2}  \omega_j\}_{j\in I}$ is an R-dual of type I of $\{S^{-1/2} f_i\}_{i\in I},$ which is a tight frame with frame bound 1; thus, Theorem \ref{207b} implies that
$\{S^{-1^/2}  \omega_j\}_{j\in I}$ is a tight Riesz sequence with bound 1, i.e., an orthonormal system.

For the proof of the other implication, if $\{ S^{-1/2} \omega_j\}_{j\in I}$ is an orthonormal system and \eqref{273a} holds, we can repeat the  proof for (ii) with the operator $S_\Omega$ replaced by $S,$ and arrive at
\begin{equation*}
S^{-1/2}  \omega_j= \sumgri \< f_i, S^{-1/2}e_j\>  h_i, \ j\in I.
\end{equation*} Applying the operator $S^{1/2}$ to this proves (ii). Finally, (iii) is a direct consequence of  Lemma \ref{213b} and (ii).\ep

Theorem \ref{213a} has several immediate consequences.  First,
we can now prove the claimed result for Gabor frames:

\bc \label{304d}
Let $\mts$ be a Gabor frame for  $L^2(\mr)$.
Then \\ $\{\frac{1}{\sqrt{ab}} E_{m/a}T_{n/b} g\}_{m,n\in \mz}$ can be realized as an
R-dual of type III of $\mts$.
 \ec

\bp Let $T$ denote the synthesis operator associated with the frame $\mts.$ Since
$\{\frac{1}{\sqrt{ab}} E_{m/a}T_{n/b} g\}_{m,n\in \mz}$ is a Riesz sequence with the same frame bounds as $\mts,$ Theorem \ref{213a} says that it is enough to show that
\bee \label{223a}  {\rm dim (ker} T)= {\rm dim (span} \{\frac{1}{\sqrt{ab}} E_{m/a}T_{n/b} g\}_{m,n\in \mz}^\perp).\ene
Since $\mts$ is a frame for $\ltr,$ we have $ab\le 1.$  If $ab<1,$ then $\mts$ is
an overcomplete frame, which according to Lemma \ref{293a} has infinite excess, i.e., $ {\rm dim (ker} \, T)= \infty.$
Also, $\{\frac{1}{\sqrt{ab}} E_{m/a}T_{n/b} g\}_{m,n\in \mz}$ is a Riesz sequence with infinite deficit, which shows that \eqref{223a} holds. On the other hand, if
$ab=1,$ then $\mts$ is a Riesz basis, so $ {\rm dim (ker} T)= 0.$ Also, for $ab=1,$
$\mts= \{\frac{1}{\sqrt{ab}} E_{m/a}T_{n/b} g\}_{m,n\in \mz},$ so clearly
${\rm dim (span} \{\frac{1}{\sqrt{ab}} E_{m/a}T_{n/b} g\}_{m,n\in \mz}^\perp)=0.$ Thus
\eqref{223a} holds for any Gabor frame $\mts$ and the proof is completed. \ep

With the insight gained by Theorem \ref{213a} we can now provide further relations between the R-duals of type I, II, and III. First we give an easy example of a frame
$\{f_i\}_{i=1}^\infty$ with an R-dual of type I which is not an
R-dual of type II, as well as a frame with an R-dual of type III which is neither an R-dual of type I nor an R-dual of type II.

\bex \label{293b}
(i) Consider the frame $\{f_i\}_{i=1}^\infty$ in Example \ref{exb1} and let $\{\nu_j\}_{j=1}^\infty$ be
the R-dual of type I of $\{f_i\}_{i=1}^\infty$ with respect to $\{e_i\}_{i=1}^\infty$ and $\{h_i\}_{i=1}^\infty=\{e_i\}_{i=1}^\infty. $
Then $\{\nu_j\}_{j=1}^\infty =\{e_1+e_2, e_3, e_4, e_5, \ldots\}.$ A calculation shows that
$\{S^{-1/2}\nu_j\}_{j=1}^\infty  =\{\frac{1}{\sqrt{2}}e_1+e_2, e_3, e_4, e_5, \ldots\},$
which is not orthonormal. Thus,
Theorem \ref{213a}(i) implies that  $\{\nu_j\}_{j=1}^\infty$ is not an R-dual
of
type II of $\{f_i\}_{i=1}^\infty.$

\noindent(ii)
Consider the frame $\{f_i\}_{i=1}^\infty$ and the Riesz basis $\{g_j\}_{j=1}^\infty$ in Example \ref{293g}. Denote the frame operator of $\{f_i\}_{i=1}^\infty$  by $S$. 
Since $S^{-1/2}g_3=(1/\sqrt{2})z_3$, the sequence $\{S^{-1/2}g_j\}_{j=1}^\infty$ is not orthonormal, which by
Theorem \ref{213a}(i) implies that $\{g_i\}_{i=1}^\infty$  is not an R-dual of type II of $\{f_i\}_{i=1}^\infty.$
By Theorem \ref{213a}(ii),  $\{g_j\}_{j=1}^\infty$ is an R-dual of type III of $\{f_i\}_{i=1}^\infty.$
\ep \enx

If $\otj$ is an R-dual of type III of a frame $\fti$ with respect to orthonormal bases
$\eti, \hti$ and a bounded bijective operator $Q,$ there is a natural way to define an R-dual of the canonical dual frame $\{\si f_i\}_{i\in I}.$ In fact, the frame operator associated with $\{\si f_i\}_{i\in I}$ is $\si,$ and the bijective operator $\widetilde{Q}:=(Q^*)^{-1}$ satisfies that
\bes || \widetilde{Q}|| = || Q^{-1}|| \le \sqrt{|| \si ||}, \ \mbox{and} \, \,
|| \widetilde{Q}^{-1}|| = || Q|| \le \sqrt{|| (\si)^{-1} ||}.\ens We have already noticed that $\{ (Q^*)^{-1} h_i\}_{i\in I}$ is the canonical dual Riesz basis of $\hti.$
We call the R-dual of type III of $\{\si f_i\}_{i\in I}$ with respect to the
orthonormal bases $\eti, \hti$ and the operator $(Q^*)^{-1}$ for the {\it canonical R-dual}
of  type III of $\{\si f_i\}_{i\in I}.$ Specifically, it is the sequence $\{\gamma_j\}_{j\in I},$ where
\bee \label{293d} \gamma_j= \sum_{i\in I} \la \si f_i, S^{1/2}e_j\ra (Q^*)^{-1}h_i=
\sum_{i\in I} \la S^{-1/2} f_i, e_j\ra (Q^*)^{-1}h_i
.\ene

For R-duals of type III we will now prove an analogue of
Proposition \ref{generalrdual}.

\bt \label{273b}

Let $\seqgri[f]$ be a frame   and  $\seqgri[\omega]$ an R-dual of $\fti$ of type III,  w.r.t. orthonormal bases
$\eti, \hti$ and a bounded bijective operator $Q.$ Denote the frame operator of $\fti$ by $S.$ Then the following hold:

\bei
\item[(i)]
$ f_i =  \sumgrj \<\omega_j, (Q^*)^{-1} h_i\> S^{1/2}e_j, \ \forall i\in I.$

\item[(ii)] The R-dual of type III of $\seqgri[f]$ with respect to some orthonormal bases
$\seqgri[e]$, $\seqgri[h]$, and an operator $Q$, is biorthogonal to the
canonical R-dual of type III of $\{\si f_i\}_{i\in I}.$
\eni \et
\bp (i) Using that $\la f_i, S^{-1/2}e_j\ra= \la S^{-1/2}f_i, e_j\ra= \la \omega_j, (Q^*)^{-1}e_j\ra,$ we have
$f_i = \sum_{j\in I} \la f_i, S^{-1/2}e_j\ra S^{1/2}e_j
=  \sum_{j\in I}  \la \omega_j, (Q^*)^{-1}e_j\ra    S^{1/2}e_j,$ as claimed.

\noindent(ii) Using the expression for  the canonical R-dual of type III in
\eqref{293d},
\begin{eqnarray*}
\<\omega_j,\gamma_n\>
&=& \<\sumgri \<  f_i, S^{-1/2} e_j\> Q h_i, \sum_{k\in I} \<  \si f_k, S^{1/2} e_n\> (Q^*)^{-1} h_k\>\\
&=& \sumgri  \sum_{k\in I}  \<  f_i, S^{-1/2} e_j\> \<  S^{1/2} e_n, \si f_k\>  \<Q h_i,(Q*)^{-1}h_k\>\\
&=& \sumgri  \<  f_i, S^{-1/2} e_j\> \<  S^{1/2} e_n, \si f_i\> \\
&= & \<  \sumgri  \<  S^{1/2} e_n, \si {f}_i\>  f_i, S^{-1/2} e_j\>
=  \<S^{1/2} e_n, S^{-1/2} e_j\> = \delta_{n,j},
\end{eqnarray*} as desired. \ep

Observe that Theorem \ref{273b} (iv) only claims biorthogonality between a given R-dual of type III of $\fti$ and the canonical R-dual of type III of $\{\si f_i \}_{i=1}^\infty.$ That is, in contrast with the situation for R-duals of type I, see Theorem \ref{207b}, we do not deal with arbitrary dual frames. The following simple example shows that the biorthogonality actually might fail in that case:

\begin{ex} \label{exb2} We return to Example \ref{exb1}, where we considered an
orthonormal basis
$\{e_i\}_{i=1}^\infty$  and the frame
$\{f_i\}_{i=1}^\infty=\{e_1, e_1, e_2, e_3, e_4, \ldots\}$. The dual frames are
exactly the frames on the form
$\{g_i\}_{i=1}^\infty= \{\alpha e_1, (1-\alpha)e_1, e_2, e_3, e_4, \cdots\},$
for some $\alpha \in \mc;$ the R-dual  of type III of $\{g_i\}_{i=1}^\infty$
with respect to $\{e_i\}_{i=1}^\infty$ and $\{h_i\}_{i=1}^\infty=\{e_i\}_{i=1}^\infty$ is
$\{\gamma_j\}_{j=1}^\infty= \{\alpha e_1+ \frac{1-\alpha}{\sqrt{|\alpha|^2+ |1-\alpha|^2}}\, e_2, e_3, e_4, \dots\}.$ In Example \ref{exb1} we found an
R-dual  of type III, $\{\omega_j\}_{j=1}^\infty,$ of $\{f_i\}_{i=1}^\infty. $ It is easy to
see that $\{\omega_j\}_{j=1}^\infty$ and $\{\gamma_j\}_{j=1}^\infty$ are biorthogonal if and only if
$\alpha +  \frac{1-\alpha}{\sqrt{2}\,\sqrt{|\alpha|^2+ |1-\alpha|^2}}=1;$
for $\alpha \in \mr,$ this holds if and only if $\alpha=1$ or $\alpha = 1/2.$ Thus, in general
$\{\omega_j\}_{j=1}^\infty$ and $\{\gamma_j\}_{j=1}^\infty$ are not biorthogonal.
\ep \enx

For frames $\fti$ and Riesz sequences $\otj$ with exactly the same bounds, we can now show that $\otj$ is an R-dual of type III of $\fti$ if and only  if $\fti$ is an R-dual of type III of $\fti.$ Again, this is a property that resembles what we know for R-duals of type I from Theorem \ref{207b}(i).

\begin{prop} \label{symm}
Let $\seqgri[f]$ be a frame  and let $\seqgrj[\omega]$ be a Riesz sequence with the same optimal bounds.
Then $\seqgrj[\omega]$ is an R-dual of type III of $\seqgri[f]$ if and only if
$\seqgri[f]$ is an R-dual of type III of $\seqgrj[\omega]$.
\end{prop}
\bp Denote the frame operators of  $\seqgri[f]$ and $\seqgrj[\omega]$ by $S$ and $S_\Omega$, respectively.
First assume that $\seqgrj[\omega]$ is an R-dual of type III of $\seqgri[f]$ with respect to some triplet ($\seqgri[e]$, $\seqgri[h]$, $Q$).
It follows from the proof of Theorem \ref{213a}(ii) that $\seqgrj[\omega]$ is an R-dual of type III of $\seqgri[f]$ with respect to the triplet ($\seqgri[e]$, $\seqgri[h]$, $\widetilde{S_\Omega^{1/2}}$), where $\widetilde{S_\Omega^{1/2}}$ is the operator defined by Lemma \ref{203a}.
Now Theorem \ref{273b}(i)
and Lemma \ref{203a} imply that
\bes f_i =  \sumgrj \< \omega_j, (\widetilde{S_\Omega^{1/2}})^{\ -1} h_i\> S^{1/2}e_j
= \sumgrj \< S_\Omega^{-1/2}\omega_j,  h_i\> S^{1/2}e_j, \ \forall i\in I.\ens
Furthermore,
$\|S^{1/2}\|=\sqrt{\|S\|}=\sqrt{\|S_\Omega\|} \ \mbox{ and } \ \|S^{-1/2}\| = \sqrt{\|S^{-1}\|}=\sqrt{\|S_\Omega^{-1}\|},$
which implies that $\seqgri[f]$ is an R-dual of type III of $\seqgrj[\omega]$ with respect to
the triplet ($\seqgri[h]$, $\seqgri[e]$, $S^{1/2}$).

Now assume that $\seqgri[f]$  is an R-dual of type III of $\seqgrj[\omega]$. 
Using techniques as in the proof of Theorem \ref{213a}(ii), one can prove that  (\ref{273a}) holds. Now Theorem \ref{213a}(ii) implies that   $\seqgrj[\omega]$  is an R-dual of type III of $\seqgri[f]$.
\ep

 Notice that when $\seqgri[f]$ is a frame for $\h$  and  $\seqgrj[\omega]$ is a Riesz sequence in $\h$ with the same optimal bounds, one has a "symmetry" representation: if (\ref{273a}) holds, then there exist
orthonormal bases $\seqgri[e]$ and $\seqgri[h]$ for $\h$ so that
$\seqgrj[\omega]$  is the $R$-dual of type III of $\seqgri[f]$ with respect
to the triplet ($\seqgri[e]$, $\seqgri[h]$, $\widetilde{S_\Omega^{1/2}}$),
and
$\seqgri[f]$ is the R-dual of type III of $\seqgrj[\omega]$ with respect to
the triplet ($\seqgri[h]$, $\seqgri[e]$, $S^{1/2}$), where $S$ and $\widetilde{S_\Omega^{1/2}}$ are as in Proposition \ref{symm}.

\vspace{.1in}
\noindent{\bf Acknowledgment:} Diana Stoeva  was supported by the 
 Austrian Science Fund (FWF) START-project FLAME ('Frames and Linear Operators for Acoustical Modeling and Parameter Estimation'; Y 551-N13).
She is grateful for the hospitality of the Technical University of Denmark, where most of the work on the paper was done.

\begin{tabbing}
text-text-text-text-text-text-text-text-text-text \= text \kill \\
Diana T. Stoeva  \> Ole Christensen \\
Acoustics Research Institute \> Technical University of Denmark\\
Wohllebengasse 12-14 \> DTU Compute\\
 Vienna A-1040, Austria\> 2800 Lyngby, Denmark  \\

Email: dstoeva@kfs.oeaw.ac.at \> Email: ochr@dtu.dk
\end{tabbing}


\begin{thebibliography}{10}

\frenchspacing

\bibitem{BCH} Balan, R., Casazza, P., and Heil, C.: {\it Deficits and excesses of frames.} Adv. Comp. Math. {\bf 18} (2003), 93--116.

\bibitem{CKL} Casazza, P., Kutyniok, G., and Lammers, M.: {\it
Duality principles in abstract frame theory.} J. Fourier Anal.
Appl. {\bf 10} no. 4, (2004), 383--408.

\bibitem{CBN} Christensen, O.:
{\it Frames and bases. An introductory course.} Birkh\"auser 2008.

\bibitem{CKK} Christensen, O.,  Kim, H.O., and Kim, R.Y.: {\it On the duality principle by Casazza,
Kutyniok, and Lammers.} J. Fourier Anal. Appl. {\bf 17} (2011), 640--655.


\bibitem{CXZ} Christensen, O., Xiao, X. C., and Zhu, Y. C.: {\it Characterizing R-duality in
Banach spaces.}
Acta Mathematica Sinica, Eng. Series {\bf 29 no.1} (2013), 75--84.

\bibitem{Da2} Daubechies, I.: {\it Ten lectures on wavelets.}
SIAM, Philadelphia, 1992.

\bibitem{DLL} Daubechies, I., Landau, H. J., and Landau, Z.: {\it
Gabor time-frequency lattices and the Wexler-Raz identity.} J.
Fourier Anal. Appl. {\bf 1} (1995), 437--478.

\bibitem{H} Dutkay, D., Han, D., and Larson, D.: {\it A duality principle for groups.}
J. Funct. Anal. {\bf 257} (2009), 1133--1143.

\bibitem{FS} Fan, Z., and Shen, Z.: {\it Dual Gramian analysis: duality principle and
unitary extension principle.} Preprint, 2013.

\bibitem{G2} Gr\"{o}chenig, K.: {\it Foundations of
time-frequency analysis.} Birkh\"{a}user, Boston, 2000.


\bibitem{Jan5} Janssen, A. J. E. M.: {\it Duality and biorthogonality
for Weyl-Heisenberg frames.} J. Fourier Anal. Appl. {\bf 1} no. 4
(1995), 403--436.




\bibitem{RoSh1} Ron, A. and Shen, Z.: {\it Weyl-Heisenberg systems and Riesz
bases in $L^2(\mr^d)$.} Duke Math. J. {\bf 89} (1997), 237--282.

\bibitem{XZ} Xiao, X. M. and Zhu, Y. C.: {\it Duality principles of frames in Banach spaces.} Acta. Math. Sci. Ser. A. Chin. {\bf 29} (2009), 94--102


\bibitem{Y} Young, R.: {\it An introduction to nonharmonic
Fourier
series.} Academic Press, New York, 1980 (revised first edition 2001).


\end{thebibliography}
\end{document}